\documentclass[a4paper,10pt,reqno]{amsart}

\usepackage[numbers]{natbib}
\usepackage{amsmath,amssymb}
\usepackage[foot]{amsaddr}
\usepackage{caption}
\usepackage{mathtools}
\usepackage[latin1]{inputenc}
\usepackage{mathrsfs}
\usepackage{enumerate}
\usepackage{textcmds}
\usepackage{graphicx,epstopdf}
\usepackage[caption=false]{subfig}
\usepackage[pdfusetitle]{hyperref}
\usepackage[capitalize]{cleveref}
\usepackage{booktabs}

\newtheorem{Theorem}{Theorem}
\newtheorem{Lemma}[Theorem]{Lemma}

\newenvironment{Proof}{\begin{proof}}{\end{proof}}

\newcommand{\D}[2]{\frac{\mathrm{d}#1}{\mathrm{d}#2}}

\newcommand{\dz}{\,\mathrm{d}z}

\newcommand{\A}{^\ast}
\newcommand{\T}{^{\operatorname{T}}}
\newcommand{\pD}[2]{\frac{\partial#1}{\partial#2}}
\newcommand{\Osym}{\text{\usefont{OMS}{cmsy}{m}{n}O}}

\newcommand{\R}{\mathbb{R}}
\newcommand{\C}{\mathbb{C}}

\newcommand{\Cz}{\C^{2}}

\newcommand{\MzC}{\mathrm{M}_{2}(\C)}
\newcommand{\LzW}{\mathrm{L}_W^2((0,1),\Cz)}
\newcommand{\ACloc}{\mathrm{AC_{loc}}((0,1),\Cz)}
\newcommand{\Tmin}{T_0}
\newcommand{\Tcls}{\overline{\Tmin}}
\newcommand{\Tmax}{T}
\newcommand{\E}{\mathrm{e}}
\newcommand{\term}[1]{\qq{#1}}
\newcommand{\U}{\mathbf{u}}

\title{On the eigenvalues of the spheroidal wave equation}
\author{Harald Schmid}
\email{h.schmid@oth-aw.de}
\address{University of Applied Sciences Amberg-Weiden, Amberg, Germany}
\keywords{spheroidal eigenvalues, Coulomb spheroidal wave equation, linear Hamiltonian system, inviscid forced Burgers' equation, deformation method}
\subjclass{33E10, 33F05, 34L15, 35F20}

\begin{document}

\begin{abstract}
This paper presents some new results on the eigenvalues of the spheroidal wave equation. We study the angular and Coulomb spheroidal wave equation as a special case of a more general linear Hamiltonian system depending on three parameters. We prove that the eigenvalues of this system satisfy a first-order quasilinear partial differential equation with respect to the parameters. This relation offers a new insight on how the eigenvalues of the spheroidal wave equation depend on the spheroidal parameter. Apart from analytical considerations, the PDE we obtain can also be used for a numerical computation of spheroidal eigenvalues.
\end{abstract}

\maketitle

\section{Introduction}
\label{sec:Introduction}

The angular spheroidal wave equation
\begin{equation} \label{ASWE}
\D{}{x}\left((1-x^2)\D{}{x}w(x)\right) + \left(\lambda + \gamma^2(1-x^2) - \frac{\mu^2}{1-x^2}\right)w(x) = 0,\quad -1<x<1
\end{equation}
appears in many fields of physics and engineering like quantum mechanics, electromagnetism, signal processing etc. Here $\mu$ is supposed to be some given real number, and since only $\mu^2$ occurs in \eqref{ASWE}, we may assume $\mu\geq 0$ without loss of generality. In many applications the so-called spheroidal parameter (or size parameter) $\gamma$ is either real or purely imaginary, so that $\gamma^2$ is a real number which may also be negative. If $\mu$ is an integer and $\gamma^2$ is real, then the separation of the Helmholtz equation in prolate ($\gamma^2>0$) or oblate ($\gamma^2<0$) spheroidal coordinates results in a second order ODE of the form \eqref{ASWE}. A slightly more general differential equation is the Coulomb spheroidal wave equation (CSWE)
\begin{equation} \label{CSWE}
\D{}{x}\left((1-x^2)\D{}{x}w(x)\right) + \left(\lambda + \beta x + \gamma^2(1-x^2) - \frac{\mu^2}{1-x^2}\right)w(x) = 0
\end{equation}
It differs from \eqref{ASWE} only by the presence of a linear term $\beta x$ with some additional parameter $\beta$, which we also assume to be real. The numbers $\lambda\in\C$ for which \eqref{CSWE} has a nontrivial bounded solution $w(x)$ on $(-1,1)$ are the eigenvalues (or characteristic values) of the CSWE, and the corresponding eigenfunctions $w(x)$ are called Coulomb spheroidal functions. They appear in astrophysics and molecular physics and provide, for example, exact wave functions for a one-electron diatomic molecule with fixed nuclei (see \cite[Chapter 9]{Falloon:2001} or \cite{KO:2018}). In our subsequent considerations, it does not make much difference whether we include the term $\beta x$ or not. Therefore, in this paper we will mainly deal with equation \eqref{CSWE}, while the results we obtain are obviously applicable to the angular spheroidal wave equation \eqref{ASWE} as well.

For technical reasons we introduce the parameters $u_1$, $u_2$, $u_3$ which are related to $\beta$, $\gamma^2$, $\lambda$ in \eqref{CSWE} by 
\begin{equation} \label{Param}
\gamma^2 = u_1,\quad\beta = -u_3-2(\mu+1)u_2,\quad\lambda = u_3+\mu(\mu+1)
\end{equation}
Moreover, we set $\alpha := \frac{1}{2}(\mu+1)$. With these parameters in mind, we will see that the $2\times 2$ differential system
\begin{equation} \label{GSWE}
\begin{pmatrix*}[r] 0 & -1 \\[1ex] 1 & 0 \end{pmatrix*}y'(z) 
- \begin{pmatrix} 2 & 2 u_2 - \frac{\alpha}{z} + \frac{\alpha}{1-z} \\[1ex] 
2 u_2 - \frac{\alpha}{z} + \frac{\alpha}{1-z} & 2(u_1+u_2^2) + \frac{u_3}{z} \end{pmatrix}
y(z) = \Lambda\begin{pmatrix} \frac{1}{1-z} & 0 \\[1ex] 0 & \frac{1}{z} \end{pmatrix}y(z),\quad 0<z<1
\end{equation}
is related to the Coulomb spheroidal wave equation by the following means: In section 2 we show that \eqref{GSWE}, in combination with appropriate boundary conditions, can be written as an eigenvalue problem for a self-adjoint differential operator $T=T(u_1,u_2,u_3)$. The eigenvalues $\Lambda$ of $T$ depend analytically on the parameters $(u_1,u_2,u_3)$, and in addition they have the following properties:
\begin{enumerate}[(a)]
\item $\lambda = u_3+\mu(\mu+1)$ is an eigenvalue of the CSWE \eqref{CSWE} for the parameters $\gamma^2 = u_1$, $\beta = -u_3-2(\mu+1)u_2$ if and only if $\Lambda=0$ is an eigenvalue of the linear Hamiltonian system \eqref{GSWE}.
\item The eigenvalues $\Lambda=\Lambda(u_1,u_2,u_3)$ of \eqref{GSWE} satisfy the first-order quasilinear partial differential equation
\begin{equation} \label{QPDE}
\begin{split}
2 u_1\frac{\partial\Lambda}{\partial u_1} 
& + \left((\Lambda+2)(u_1+u_2^2)+\Lambda+u_2+u_3\right)\frac{\partial\Lambda}{\partial u_2} \\[-1ex]
& + \left((2\Lambda+u_3+2)(2\Lambda u_2-2\mu-1)-2(\mu+1)(u_1+u_2^2)+2\mu\right)\frac{\partial\Lambda}{\partial u_3} \\
& = (1+2\mu-2\Lambda u_2)(\Lambda+2)-2\mu
\end{split}
\end{equation}
\end{enumerate}
According to (a), the eigenvalue problem \eqref{GSWE} may be regarded as a generalization of the Coulomb spheroidal wave equation, and the PDE \eqref{QPDE} describes an analytic relation between the eigenvalues and the parameters, where the eigenfunctions need not to be known. A similar observation has already been made in \cite{BSW:2005} for the Chandrasekhar-Page equation, and \cite{Schmid:2021} provides a method by which one can establish such a relation between the eigenvalues and the parameters in terms of a PDE for other linear Hamiltonian systems as well. In the present paper we will use this \term{deformation method} from \cite{Schmid:2021} to prove assertion (b). Finally, in section 3 we solve \eqref{QPDE} by the method of characteristics, where we consider the prolate spheroidal wave equation ($\beta=0$, $\gamma^2>0$) in more detail. In this case, we can associate a linear $2\times 2$ differential system to \eqref{ASWE} which contains only two parameters in addition to the eigenvalue parameter, so that the partial differential equation for the eigenvalues can be further simplified.

\section{A linear Hamiltonian system associated to the CSWE}
\label{sec:Investigation}

In this section we study the linear $2\times 2$ Hamiltonian system \eqref{GSWE} on the interval $(0,1)$. It depends on three parameters $u_1,u_2,u_3\in\R$, where $\alpha=\frac{1}{2}(\mu+1)\geq\frac{1}{2}$ with some constant $\mu\geq 0$, and the prime denotes the derivative with respect to $z$. If we define
\begin{equation*}
J := \begin{pmatrix*}[r] 0 & -1 \\[1ex] 1 & 0 \end{pmatrix*},\quad
H(z,\U) := \begin{pmatrix} 2 & -\frac{\alpha}{z}+\frac{\alpha}{1-z}+2u_2 \\[1ex] 
-\frac{\alpha}{z}+\frac{\alpha}{1-z}+2u_2 & \frac{u_3}{z}+2(u_1+u_2^2) \end{pmatrix},\quad
W(z) := \begin{pmatrix} \frac{1}{1-z} & 0 \\[1ex] 0 & \frac{1}{z} \end{pmatrix}
\end{equation*}
and
\begin{equation} \label{tau}
\tau y := W(z)^{-1}\left(J y'(z) - H(z,\U)y(z)\right)
\end{equation}
where $\mathbf{u} := (u_1,u_2,u_3)\in\R^3$ combines the three parameters, then \eqref{GSWE} can be written in the form $\tau y = \Lambda y$. Moreover, as $W(z)\A=W(z)>0$ and $H(z,\U)\A=H(z,\U)$ for all $z\in(0,1)$, $\tau$ defines for fixed parameters $\U\in\R^3$ a formally self-adjoint differential expression in the Hilbert space
\begin{equation*}
\LzW := \Big\{f:(0,1)\longrightarrow\Cz\ \big|\ \int_0^1 f(z)\A W(z) f(z)\dz < \infty\Big\}
\end{equation*}
of square-integrable vector functions to the weight function $W(z)$ with scalar product
\begin{equation*}
\langle f,g\rangle_W := \int_0^1 f(z)\A W(z)\,g(z)\dz
\end{equation*}
The maximal operator $\Tmax y := \tau y$ generated by $\tau$ has the domain (cf. \cite[Section 3]{Weidmann:1987})
\begin{equation} \label{Domain}
D(\Tmax) := \left\{y\in\LzW\ :\ y\in\ACloc\mbox{ and }\tau y\in\LzW\right\}
\end{equation}
whereas the domain of the minimal operator $\Tmin y := \tau y$ associated to $\tau$ is given by
\begin{equation*}
D(\Tmin) := \big\{y\in D(\Tmax)\,:\,y\mbox{ has compact support in }(0,1)\big\}
\end{equation*}
In order to establish the self-adjointness of $T$, we first prove that $\tau$ is in the limit point case at $z=0$ and $z=1$. To this end, we consider the differential equation $\tau y = \Lambda y$ in the complex plane. It is equivalent to the $2\times 2$ system
\begin{equation} \label{SysA}
y'(z) = \left(
\frac{1}{z}\begin{pmatrix} -\alpha & u_3+\Lambda \\[1ex] 0 & \alpha \end{pmatrix} 
+ \frac{1}{1-z}\begin{pmatrix} \alpha & 0 \\[1ex] -\Lambda & -\alpha \end{pmatrix}
+ \begin{pmatrix} 2 u_2 & 2(u_1+u_2^2) \\[1ex] -2 & - 2 u_2 \end{pmatrix}\right)y(z)
\end{equation}
with regular singular points at $z=0$ and $z=1$. If $\mathfrak{D}_0:=\{z\in\C:|z|<1\}$ denotes the unit disk  centered at $0$, then \cite[Lemma 6]{BSW:2005} implies that for fixed $(\Lambda,\U)\in\C\times\R^3$ there exists a fundamental matrix of the form $Y_0(z) = F_0(z,\Lambda,\U)z^A z^{Q_0(\Lambda,\U)}$, $z\in\mathfrak{D}_0$, where $F_0:\mathfrak{D}_0\times\C\times\R^3\longrightarrow\MzC$ is an analytic matrix function,
\begin{equation*}
F_0(0,\Lambda,\U) = \begin{pmatrix} 1 & \frac{u_3+\Lambda}{2\alpha} \\[1ex] 0 & 1 \end{pmatrix},\quad
A := \begin{pmatrix} -\alpha & 0 \\[1ex] 0 & \alpha \end{pmatrix},\quad
Q_0(\Lambda,\U) = \begin{pmatrix} 0 & 0 \\[1ex] q_0(\Lambda,\U) & 0 \end{pmatrix}
\end{equation*}
and $q_0:\C\times\R^3\longrightarrow\C$ is an analytic function with $q_0\equiv 0$ if $2\alpha=\mu+1$ is not an integer. Furthermore, the transformation $\tilde y(z) = y(1-z)$ yields
\begin{equation*}
\tilde y'(z) = \left(
\frac{1}{z}\begin{pmatrix} -\alpha & 0 \\[1ex] \Lambda & \alpha \end{pmatrix} 
+ \frac{1}{1-z}\begin{pmatrix} \alpha & -u_3-\Lambda \\[1ex] 0 & -\alpha \end{pmatrix}
+ \begin{pmatrix} -2u_2 & -2(u_1+u_2^2) \\[1ex] 2 & 2u_2 \end{pmatrix}\right)\tilde y(z)
\end{equation*}
and hence \eqref{SysA} also possesses a fundamental matrix of the form
$Y_1(z) = F_1(z,\Lambda,\U)(1-z)^A (1-z)^{Q_1(\Lambda,\U)}$, $z\in\mathfrak{D}_1:=\{z\in\C:|z-1|<1\}$,
where $F_1:\mathfrak{D}_1\times\C\times\R^3\longrightarrow\MzC$ is an analytic matrix function,
\begin{equation*}
F_1(1,\Lambda,\U) = \begin{pmatrix} 1 & 0 \\[1ex] -\frac{\Lambda}{2\alpha} & 1\end{pmatrix},\quad
Q_1(\Lambda,\U) = \begin{pmatrix} 0 & 0 \\[1ex] q_1(\Lambda,\U) & 0 \end{pmatrix}
\end{equation*}
and $q_1:\C\times\R^3\longrightarrow\C$ being an analytic function satisfying $q_1\equiv 0$ if $2\alpha=\mu+1$ is not an integer. For fixed parameters $(\Lambda,\U)\in\C\times\R^3$, the solution
\begin{equation*}
y_0(z) := Y_0(z)e_1 
= F_0(z)\begin{pmatrix} z^{-\alpha} & 0 \\[1ex] 0 & z^\alpha \end{pmatrix}
\begin{pmatrix} 1 & 0 \\[1ex] q_0\log z & 1 \end{pmatrix}\begin{pmatrix} 1 \\[1ex] 0 \end{pmatrix}
\end{equation*}
of \eqref{SysA} asymptotically behaves like $y_0(z)\sim z^{-\alpha}e_1$ as $z\to 0$ (here and in the following, $\{e_1,e_2\}$ denotes the standard basis of $\Cz$). Since $2\alpha\geq 1$ and $y_0(z)\A W(z) y_0(z)\sim z^{-2\alpha}$ as $z\to 0$, the function $y_0(z)\A W(z) y_0(z)$ is not integrable at $z=0$, and hence $y_0$ does not lie left in $\LzW$. Moreover, for the solution
\begin{equation*}
y_1(z) := Y_1(z)e_1 
= F_1(z)\begin{pmatrix} (1-z)^{-\alpha} & 0 \\[1ex] 0 & (1-z)^\alpha \end{pmatrix}
\begin{pmatrix} 1 & 0 \\[1ex] q_1\log(1-z) & 1 \end{pmatrix}\begin{pmatrix} 1 \\[1ex] 0 \end{pmatrix}
\end{equation*}
we get $y_1(z)\A W(z) y_1(z)\sim (1-z)^{-2\alpha-1}$ as $z\to 1$. Therefore, $y_1(z)\A W(z) y_1(z)$ is not integrable at $z=1$, and thus $y_1$ does not lie right in $\LzW$. According to Weyl's alternative (see e.\,g. \cite[Theorem 5.6]{Weidmann:1987}), the differential expression $\tau$ is in the limit point case at both endpoints. Hence, the closure of its associated minimal operator $\Tcls$ is a self-adjoint extension of $T_0$ (cf. \cite[Theorem 5.4]{SunShi:2010}), and because of $\Tcls = \Tcls\A = \Tmax$ (see \cite[Theorem 3.9]{Weidmann:1987}), this extension coincides with the maximal operator $T$, so that $T$ is actually a self-adjoint differential operator in the Hilbert space $\LzW$. 

In the following, we will study the eigenvalues of the self-adjoint operator $T=T(\U)$ associated to $\tau$ and their dependence on the parameters $\U\in\R^3$. For this purpose, we introduce the functions $\eta_0(z,\Lambda,\U) := Y_0(z,\Lambda,\U)e_2$ and $\eta_1(z,\Lambda,\U) := Y_1(z,\Lambda,\U)e_2$, which asymptotically behave like
\begin{equation} \label{Asym}
\eta_0(z,\Lambda,\U)\sim z^\alpha\begin{pmatrix} \frac{u_3+\Lambda}{2\alpha} \\[1ex] 1 \end{pmatrix}
\quad\mbox{as $z\to 0$},\quad
\eta_1(z,\Lambda,\U)\sim (1-z)^\alpha\begin{pmatrix} 0 \\[1ex] 1 \end{pmatrix}\quad\mbox{as $z\to 1$}
\end{equation}
Note that $y_0$, $\eta_0$ as well as $y_1$, $\eta_1$ form a fundamental system of \eqref{GSWE}. Since the solutions $y_0(z) = Y_0(z)e_1$ and $y_1(z) = Y_1(z)e_1$ are not square-integrable with respect to $W(z)$, $\Lambda\in\R$ is an eigenvalue of $T(\U)$ for fixed $\U\in\R^3$ if and only if there exists a nontrivial solution $y$ of \eqref{GSWE} which is a constant multiple of both $\eta_0$ and $\eta_1$. Additionally, by virtue of $\alpha\geq\frac{1}{2}$ and the asymptotic behavior of the fundamental solutions at the boundary points, a nontrivial solution $y$ of \eqref{GSWE} is an eigenfunction if and only if $(z(1-z))^{-1/2}y(z)$ is bounded on $(0,1)$. Furthermore, we can write $\eta_0$ as a linear combination
\begin{equation} \label{Connect}
\eta_0(z,\Lambda,\U) = \Theta(\Lambda,\U) y_1(z,\Lambda,\U) + \Omega(\Lambda,\U)\eta_1(z,\Lambda,\U)
\end{equation}
with $\Theta$, $\Omega$ called \emph{connection coefficients}. Since for any fixed $z_0\in(0,1)$
\begin{equation*}
\begin{pmatrix} \Theta(\Lambda,\U) \\[0.5ex] \Omega(\Lambda,\U) \end{pmatrix} = Y_1(z_0,\Lambda,\U)^{-1}Y_0(z_0,\Lambda,\U)e_2
\end{equation*}
where the fundamental matrices depend analytically on the parameters, also $\Theta=\Theta(\Lambda,\U)$ and $\Omega=\Omega(\Lambda,\U)$ depend analytically on $(\Lambda,\U)\in\C\times\R^3$. Due to our considerations above, $\Lambda\in\R$ is an eigenvalue of $T(\U)$ if and only if $\Theta(\Lambda,\U)=0$. As $\Theta(\Lambda,\U)$ is an entire function for fixed $\U\in\R^3$, the zeros of $\Theta(\Lambda,\U)$ are isolated, and thus also the eigenvalues of $T(\U)$ are isolated with multiplicity $1$.

In a first step, we will figure out how the eigenvalues of \eqref{CSWE} and \eqref{GSWE} are related.

\begin{Lemma} \label{Embed} 
$\Lambda=0$ is an eigenvalue of the self-adjoint operator $T(u_1,u_2,u_3)$ associated to the linear Hamiltonian system \eqref{GSWE} if and only if $\lambda = u_3+\mu(\mu+1)$ is an eigenvalue of the Coulomb spheroidal wave equation \eqref{CSWE} for the parameters $\gamma^2 = u_1$ and $\beta = -u_3-2(\mu+1)u_2$.
\end{Lemma}

\begin{Proof}
The system \eqref{SysA} with $\Lambda=0$ is equivalent to
\begin{equation*}
y'(z) = \left(
\frac{1}{z}\begin{pmatrix} -\alpha & u_3 \\[1ex] 0 & \alpha \end{pmatrix} 
+ \frac{1}{1-z}\begin{pmatrix} \alpha & 0 \\[1ex] 0 & -\alpha \end{pmatrix}
+ \begin{pmatrix} 2 u_2 & 2(u_1+u_2^2) \\[1ex] -2 & - 2 u_2 \end{pmatrix}\right)y(z)
\end{equation*}
If we introduce the vector function
\begin{equation} \label{Trafo}
\hat y(x) = \begin{pmatrix} v(x) \\[0.5ex] w(x) \end{pmatrix} := 2(1-x^2)^{-\frac{1}{2}}y\left(\tfrac{x+1}{2}\right),\quad x\in(-1,1)
\end{equation}
then this differential system takes the form
\begin{equation} 
\hat y'(x) = \left(
\frac{1}{1+x}\begin{pmatrix} -\alpha-\frac{1}{2} & u_3 \\[1ex] 0 & \alpha-\frac{1}{2} \end{pmatrix} 
+ \frac{1}{1-x}\begin{pmatrix} \alpha+\frac{1}{2} & 0 \\[1ex] 0 & -\alpha+\frac{1}{2} \end{pmatrix}
+ \begin{pmatrix} u_2 & u_1+u_2^2 \\[1ex] -1 & -u_2 \end{pmatrix}\right)\hat y(x) \label{LHS0}
\end{equation}
Since $\hat y(x) = (z(1-z))^{-1/2}y(z)$ with $z=\frac{1}{2}(x+1)\in(0,1)$, we obtain that $\Lambda=0$ is an eigenvalue of $T(\U)$ if and only if \eqref{LHS0} has a nontrivial solution $\hat y(x)$ which is bounded on $(-1,1)$. Using $\alpha=\frac{1}{2}(\mu+1)$, we can write \eqref{LHS0} in terms of its components $v(x)$ and $w(x)$:
\begin{align*}
v'(x) & = \left(u_2+\frac{(\mu+2)x}{1-x^2}\right)v(x) + \left(u_1+u_2^2+\frac{u_3}{1+x}\right)w(x) \\
w'(x) & = -v(x) - \left(u_2+\frac{\mu x}{1-x^2}\right)\,w(x)
\end{align*}
Solving the second equation for $v(x)$ yields
\begin{equation} \label{SolV}
v(x) = -w'(x) - \left(u_2+\frac{\mu x}{1-x^2}\right)\,w(x)
\end{equation}
and inserting this expression into the first equation gives the differential equation
\begin{equation*}
\D{}{x}\left((1-x^2)\D{}{x}w(x)\right) + \left(u_3+\mu(\mu+1) - (u_3+2(\mu+1)u_2)x + u_1(1-x^2) - \frac{\mu^2}{1-x^2}\right)w(x) = 0
\end{equation*}
for $w(x)$ which is a CSWE \eqref{CSWE} with parameters \eqref{Param}. Conversely, any bounded solution $w(x)$ of this ODE is a Coulomb spheroidal function having the form $w(x)=(1-x^2)^{\mu/2}\varphi(x)$ with some entire function $\varphi(x)$, cf. \cite[Section 3.2, Satz 1]{MS:1954}. If we define $v(x)$ by \eqref{SolV}, then $v(x)=-(1-x^2)^{\mu/2}\varphi'(x)-u_2(1-x^2)^{\mu/2}\varphi(x)$ is also bounded on $(-1,1)$, and hence the function $\hat y$ defined by \eqref{Trafo} is a bounded solution of \eqref{LHS0} which implies that $\Lambda=0$ is an eigenvalue of \eqref{GSWE}.
\end{Proof}

In the next step we examine which way the eigenvalues $\Lambda$ of \eqref{GSWE} depend on the parameters $(u_1,u_2,u_3)$.

\begin{Theorem} \label{Main} 
Let $T=T(\U)$ be the self-adjoint realization of the differential expression \eqref{tau} with domain $D(T)$ given by \eqref{Domain}, and suppose that $\Lambda_0$ is an eigenvalue of $T=T(\U_0)$ for some given parameter vector $\U_0\in\R^3$. Then there exists a domain $U\subset\R^3$ with $\U_0\in U$ and eigenvalues $\Lambda(\U)$ of $T(\U)$ with $\Lambda(\U_0)=\Lambda_0$, such that $\Lambda=\Lambda(\U)$ depends analytically on $\U=(u_1,u_2,u_3)\in U$. In addition, $\Lambda$ satisfies the first-order partial differential equation \eqref{QPDE} on $U$.
\end{Theorem}

\begin{Proof}
We first prove the analytic dependence of the eigenvalues on the parameters $\U=(u_1,u_2,u_3)$ in a neighborhood of $\U_0=:(a,b,c)$. For this, let 
\begin{equation*}
B(z,\U) := H(z,\U) - H(z,\U_0) 
= \begin{pmatrix} 0 & 2(u_2-b) \\[1ex] 2(u_2-b) & 2(u_1+u_2^2-a-b^2)+\frac{u_3-c}{z} \end{pmatrix}
\end{equation*}
If $A(\U)$ denotes the symmetric multiplication operator $A(\U)y := -W(z)^{-1}B(z,\U)y$ in $\LzW$, then $T(\U) = T(\U_0) + A(\U)$. Since
\begin{equation*}
W(z)^{-1/2}B(z,\U)W(z)^{-1/2} = 
\begin{pmatrix} 0 & 2(u_2-b)\sqrt{z(1-z)} \\[1ex] 2(u_2-b)\sqrt{z(1-z)} & 2(u_1+u_2^2-a-b^2)z+u_3-c \end{pmatrix}
\end{equation*}
we have $|W(z)^{-1/2}B(z,\U)W(z)^{-1/2}|\leq C$ for all $z\in(0,1)$ with some constant $C=C(\U)\geq 0$, and
\begin{align*}
\|A y\|_W^2 
& = \int_0^1 \left(W(z)^{-1}B(z,\U)\,y(z)\right)\A W(z)\left(W(z)^{-1}B(z,\U)\,y(z)\right)\dz \\
& = \int_0^1 |W(z)^{-1/2}B(z,\U)W(z)^{-1/2}W(z)^{1/2}y(z)|^2\dz 
\leq\int_0^1 C^2|W(z)^{1/2}y(z)|^2\dz = C^2\|y\|_W^2
\end{align*}
implies that $A=A(\U)$ is even a bounded operator on $\LzW$ for all $\U\in\C^3$. From \cite[Chap. V, \S\,4, Theorem 4.3 and Chap. VII, \S\,2, Theorem 2.6]{Kato:1995} it follows that $T(\U)$ defines a holomorphic family of self-adjoint operators, where its domain $D(T)$ is independent of $\U\in\C^3$. As shown in \cite[Chap. VII, \S\,2, Sec. 4 and \S\,3, Sec. 1 -- 2]{Kato:1995} there exists a domain $U\subset\R^3$ with $\U_0\in U$ as well as an analytical function $\Lambda=\Lambda(u_1,u_2,u_3)$ with $\Lambda(\U_0)=\Lambda_0$, such that $\Lambda(\U)$ is a simple eigenvalue of $T(\U)$ for all $\U\in U$.

Now, if $\Lambda(\U)$ is an eigenvalue of $T(\U)$, then each corresponding eigenfunction is a constant multiple of the fundamental solutions $\eta_0$ and $\eta_1$, respectively. Hence, $\eta_0(z) = \Omega(\U)\cdot\eta_1(z)$ for all $z\in(0,1)$ with some factor $\Omega(\U)\neq 0$ which depends analytically on $\U\in U$, and therefore $\eta_0(z) = z^{\alpha}(1-z)^{\alpha}h(z,\U)$, where $h:\mathfrak{D}\times U\longrightarrow\Cz$ is an analytic function on some domain $[0,1]\subset\mathfrak{D}\subset\C$ which satisfies
\begin{equation*}
h(0,\U)=\begin{pmatrix} \frac{u_3+\Lambda(\U)}{2\alpha} \\[1ex] 1 \end{pmatrix}\quad\mbox{and}\quad
h(1,\U)=\begin{pmatrix} 0 \\[1ex] \Omega(\U) \end{pmatrix}
\end{equation*}
due to \eqref{Asym}. Obviously $\eta_0\in\LzW$, and
\begin{equation*}
\nu(\U) := \langle\eta_0,\eta_0\rangle_W = \int_0^1 h(z,\U)\A 
\begin{pmatrix} z^{\mu+1}(1-z)^{\mu} & 0 \\[1ex] 0 & z^{\mu}(1-z)^{\mu+1} \end{pmatrix}h(z,\U)\dz
\end{equation*}
defines a smooth function $\nu:U\longrightarrow(0,\infty)$. Setting $\psi(z,\U) := \nu(\U)^{-1/2}h(z,\U)$, then 
\begin{equation*}
y(z,\U):=\nu(\U)^{-1/2}\eta_0(z,\Lambda(\U),\U)=z^{\alpha}(1-z)^{\alpha}\psi(z,\U)
\end{equation*}
is an eigenfunction corresponding to $\Lambda(\U)$ for which $\langle y,y\rangle_W=1$ holds for all $\U\in U$, and $y=y(z,\U)$ depends infinitely differentiable on $(z,\U)\in(0,1)\times U$. Since the matrix functions $W:(0,1)\longrightarrow\MzC$ and $H:(0,1)\times\R^3\longrightarrow\MzC$ are also infinitely differentiable with respect to all variables, the assumptions (a), (b), (c) in \cite{Schmid:2021} are fulfilled. In addition, if we introduce the matrix function
\begin{equation*}
G(z,\Lambda,\U) := \begin{pmatrix} 
2 z & 2 z u_2 + \Lambda u_2 - \mu - \frac{1}{2} \\[1ex] 
2 z u_2 + \Lambda u_2 - \mu - \frac{1}{2} & 2(z-1)(u_1+u_2^2) \end{pmatrix}
\end{equation*}
then a direct calculation gives
\begin{equation} \label{DefGSWE}
\pD{G}{z} + (\Lambda W + H)JG - GJ(\Lambda W + H) = \sum_{k=1}^{3} f_k\pD{H}{u_k} + g W
\end{equation}
where the functions $f_k=f_k(\Lambda,\U)$ and $g=g(\Lambda,\U)$ on the right hand side are defined by
\begin{align*}
f_1 & := 2 u_1 \\
f_2 & := (\Lambda+2)(u_1+u_2^2)+\Lambda+u_2+u_3 \\
f_3 & := (2\Lambda+u_3+2)(2\Lambda u_2-2\mu-1)-2(\mu+1)(u_1+u_2^2)+2\mu \\
  g & := (1+2\mu-2\Lambda u_2)(\Lambda+2)-2\mu
\end{align*}
Here, $G:(0,1)\times\C\times U\longrightarrow\MzC$ and $f_k, g:\C\times U\longrightarrow\C$ are differentiable functions, and therefore also assumption (d) in \cite{Schmid:2021} is satisfied. In addition to the \term{deformation equation} \eqref{DefGSWE}, the function
\begin{align*}
\sum_{k=1}^{3} f_k(\Lambda(\U),\U)\pD{y}{u_k} + JG(z,\Lambda(\U),\U)y = 
z^\alpha(1-z)^\alpha\left(\sum_{k=1}^{3} f_k(\Lambda(\U),\U)\pD{\psi}{u_k} + JG(z,\Lambda(\U),\U)\psi\right)
\end{align*}
lies in $\LzW$. Finally, from \cite[Theorem 2.2]{Schmid:2021} it follows that $\Lambda=\Lambda(u_1,u_2,u_3)$ solves
\begin{equation*}
\sum_{k=1}^{3} f_k(\Lambda,u_1,u_2,u_3)\pD{\Lambda}{u_k} = g(\Lambda,u_1,u_2,u_3)
\end{equation*}
on $U$, and this is exactly the first-order quasilinear PDE \eqref{QPDE}.
\end{Proof}

The eigenvalues of \eqref{GSWE} coincide with the zeros of the function $\Theta(\Lambda,\U)$, and this is one of the connection coefficients appearing in \eqref{Connect}. The question arises how this function can be calculated numerically for given parameter values $\U\in\R^3$ and $\Lambda\in\C$.
In \cite{Schmid:2023}, the Coulomb spheroidal wave equation has already been associated to a $2\times 2$ differential system, which is, however, different from the linear Hamiltonian system studied in the present paper. Nevertheless, the method suggested in \cite{Schmid:2023} for calculating the connection coefficient $\Theta$ can also be applied to the system \eqref{GSWE}. To this end, we apply the transformation $\hat y(z) = (1-z)^{\alpha-1}y(z)$, which turns \eqref{GSWE} resp. \eqref{SysA} into the differential system
\begin{equation} \label{SysC}
\hat y'(z) = \left(\frac{1}{z}\,A_0 + \frac{1}{z-1}\,A_1 + C\right)\hat y(z)
\end{equation}
where
\begin{equation*}
A_0 = \begin{pmatrix} -\alpha & u_3+\Lambda \\[1ex] 0 & \alpha \end{pmatrix},\quad
A_1 = \begin{pmatrix} -1 & 0 \\[1ex] \Lambda & 2\alpha-1 \end{pmatrix},\quad
 C  = \begin{pmatrix} 2 u_2 & 2(u_1+u_2^2) \\[1ex] -2 & - 2 u_2 \end{pmatrix}
\end{equation*}
For fixed $(\lambda,\U)\in\C\times\R^3$ it has the Floquet solution
\begin{equation*}
\hat\eta_0(z) := (1-z)^{\alpha-1}\eta_0(z) = z^\alpha h(z),\quad
h(z)=\sum_{k=0}^\infty z^k d_k,\quad
d_0 := \begin{pmatrix} \frac{u_3+\Lambda}{2\alpha} \\[1ex] 1 \end{pmatrix} 
\end{equation*}
where $h$ is holomorphic in $\mathfrak{D}_0$, and $d_0$ is an eigenvector of $A_0$ corresponding to the eigenvalue $\alpha$. Finally, once we have determined the series coefficients $d_k$ using the recurrence relation
\begin{equation} \label{SphRF}
b_k := \begin{pmatrix*} 
\frac{2(k u_2-u_3)}{k(k+\mu+1)} & \frac{2k(u_1+u_2^2)-2u_2u_3}{k(k+\mu+1)} \\[1ex] -\frac{2}{k} & -\frac{2u_2}{k} \end{pmatrix*}b_{k-1} - \begin{pmatrix*} 0 & \frac{(\mu+1)u_3}{k(k+\mu+1)} \\[1ex] 0 & \frac{\mu+1}{k} \end{pmatrix*}d_{k-1},\quad d_k := d_{k-1}+b_k
\end{equation}
for $k=1,2,3,\ldots$ starting with $b_0 := d_0$, then we can compute $\Theta$ with arbitrary precision by means of

\begin{Lemma} \label{Limit}
Let $\mu\geq 0$ and $\U=(u_1,u_2,u_3)\in\R^3$ be fixed. Then for each $\Lambda\in\C$ we have 
\begin{equation*}
\Theta(\Lambda,\U) := \lim_{k\to\infty}e_1\T d_k(\Lambda,\U)
\end{equation*}
More precisely, we obtain $e_1\T d_k = \Theta + \Osym(k^{\delta-\mu-2})$ as $k\to\infty$ for arbitrary small $\delta>0$.
\end{Lemma}

This assertion, as well as formula \eqref{SphRF}, is obtained by a similar reasoning as given in \cite[Section 3]{Schmid:2023}, and therefore the proof is omitted here. If we set $\Lambda=0$ and calculate the zeros $\U=(u_1,u_2,u_3)$ of the function $\Theta(0,\U)$, then we get the eigenvalues $\lambda = u_3+\mu(\mu+1)$ of the Coulomb spheroidal wave equation \eqref{CSWE} for the parameters $\gamma^2 = u_1$ and $\beta = -u_3-2(\mu+1)u_2$ according to Lemma \ref{Embed}. This approach, like the method presented in \cite{Schmid:2023}, provides a simple and efficient way to obtain the eigenvalues of a CSWE. However, with increasing parameter values, numerical effects such as rounding errors or digit cancellation may occur, which require some computational effort to overcome. In this context the partial differential equation \eqref{QPDE} may be helpful. It not only provides us with new insights into the analytical dependence of the eigenvalues $\Lambda=\Lambda(\U)$ on the parameters, but can also be used to determine the Coulomb spheroidal eigenvalues numerically. For this purpose, we need to examine the PDE and its characteristics in more detail.

\section{Analysis of the PDE for the eigenvalues}
\label{sec:Evaluation}

The PDE \eqref{QPDE} looks very complicated at first sight, but it has some remarkable properties that simplify the calculation of the solution by the method of characteristics. The characteristic curves related to the quasilinear PDE \eqref{QPDE} satisfy the (nonlinear) first-order differential system
\begin{equation} \label{GenChar}
\begin{split}
\D{u_1}{t} & = 2 u_1(t) \\
\D{u_2}{t} & = (\Lambda(t)+2)(u_1(t)+u_2(t)^2)+\Lambda(t)+u_2(t)+u_3(t) \\
\D{u_3}{t} & = (2\Lambda(t)+u_3(t)+2)(2\Lambda(t)u_2(t)-2\mu-1)-2(\mu+1)(u_1(t)+u_2(t)^2)+2\mu \\
\D{\Lambda}{t} & = (1+2\mu-2\Lambda(t)u_2(t))(\Lambda(t)+2)-2\mu
\end{split}
\end{equation}
We can choose the parameter for the characteristic curve in such a way that the initial conditions are given at $t=0$. If we define the scalar function
\begin{equation*}
\Psi(t) := \Lambda(t)(u_1(t)+u_2(t)^2-1) - 2(\mu+1)u_2(t) - u_3(t)
\end{equation*}
then \eqref{GenChar} and a straightforward calculation yields $\D{}{t}\Psi = \Psi$. In particular, if $\Psi(0)=0$, then $\Psi\equiv 0$, and therefore $u_3(t) = \Lambda(t)(u_1(t)+u_2(t)^2-1) - 2(\mu+1)u_2(t)$ holds for all $t$ along such a characteristic curve. If, in addition, $\Lambda(t_0)=0$ for some $t_0\in\R$, then we get $u_3(t_0)=-2(\mu+1)u_2(t_0)$. According to Lemma \ref{Embed}, $\Lambda(t_0)=0$ also implies that $\lambda = u_3(t_0)+\mu(\mu+1)$ is an eigenvalue of the Coulomb spheroidal wave equation \eqref{CSWE} for the parameters $\gamma^2 = u_1(t_0)$ and $\beta = -u_3(t_0)-2(\mu+1)u_2(t_0) = 0$, where $\beta=0$ in turn reduces the CSWE to the ordinary (angular) spheroidal wave equation. All in all, the above considerations provide an alternative method for the numerical computation of spheroidal eigenvalues:

\begin{Theorem} \label{CharEqu}
Suppose that $u_1(t)$, $u_2(t)$ and $\Lambda(t)$ are solutions of the differential system
\begin{equation} \label{RedChar}
\begin{split}
\D{u_1}{t} & = 2 u_1(t) \\
\D{u_2}{t} & = 2(\Lambda(t)+1)(u_1(t)+u_2(t)^2)-(2\mu+1)u_2(t) \\
\D{\Lambda}{t} & = (1+2\mu-2\Lambda(t)u_2(t))(\Lambda(t)+2)-2\mu
\end{split}
\end{equation}
on an interval $0\in I\subset\R$ satisfying $u_1(0) = a$, $u_2(0)=b$, $\Lambda(0)=\Lambda_0$ with some given initial values $a,b\in\R$ such that $\Theta(\Lambda_0,a,b,(a+b^2-1)\Lambda_0-2(\mu+1)b) = 0$. If $\Lambda(t_0)=0$ holds for some $t_0\in I$, then $\lambda = (\mu+1)(\mu-2 u_2(t_0))$ is an eigenvalue of the angular spheroidal wave equation \eqref{ASWE} for the parameter $\gamma^2 = u_1(t_0)$.
\end{Theorem}

\begin{Proof}
If we set $c := (a+b^2-1)\Lambda_0-2(\mu+1)b$, then $\Theta(\Lambda_0,a,b,c) = 0$, and hence $\Lambda_0$ is an isolated eigenvalue of the operator $T(\U_0)$ for $\U_0:=(a,b,c)$. According to Theorem \ref{Main}, we can find an analytic continuation $\Lambda_0:U\longrightarrow\R$ on some domain $U\subset\R^3$ with $\U_0\in U$, such that $\Lambda_0(\U)$ is an eigenvalue of $T(\U)$ satisfying, in addition, the partial differential equation \eqref{QPDE}. Defining $u_3(t) := \Lambda(t)(u_1(t)+u_2(t)^2-1) - 2(\mu+1)u_2(t)$, then we can write the second equation in \eqref{RedChar} as
\begin{equation*}
\D{u_2}{t} = (\Lambda(t)+2)(u_1(t)+u_2(t)^2)+u_2(t)+u_3(t)+\Lambda(t)
\end{equation*}
Moreover, the differential equations \eqref{RedChar} imply
\begin{align*}
\D{u_3}{t} 
& = \D{\Lambda}{t}(u_1(t)+u_2(t)^2-1)+\Lambda(t)\left(\D{u_1}{t}+2u_2(t)\D{u_2}{t}\right) - 2(\mu+1)\D{u_2}{t} \\
& = (2\Lambda(t)+u_3(t)+2)(2\Lambda(t)u_2(t)-2\mu-1)-2(\mu+1)(u_1(t)+u_2(t)^2)+2\mu
\end{align*}
Hence, the functions $u_1(t)$, $u_2(t)$, $u_3(t)$, $\Lambda(t)$ satisfy the characteristic equations \eqref{GenChar}, and from the initial values $u_1(0) = a$, $u_2(0)=b$, $u_3(0)=c$, $\Lambda(0)=\Lambda_0=\Lambda_0(a,b,c)$ it follows that this characteristic curve lies completely on the integral surface $\Lambda_0(\U)$. In particular, $\Lambda(t) = \Lambda_0(u_1(t),u_2(t),u_3(t))$, and therefore $\Lambda(t)$ is an eigenvalue of \eqref{GSWE} for the parameter values $u_1(t)$, $u_2(t)$, $u_3(t)$. Finally, if $\Lambda(t_0)=0$, then $u_3(t_0) = -2(\mu+1)u_2(t_0)$ implies $\beta = -u_3(t_0)-2(\mu+1)u_2(t_0) = 0$, and $\lambda = (\mu+1)(\mu-2 u_2(t_0))$ is an eigenvalue of \eqref{ASWE} for the spheroidal parameter $\gamma^2 = u_1(t_0)$.
\end{Proof}

We may, for example, choose $b=0$ in Theorem \ref{CharEqu}, and then we need to find a zero $\Lambda_0$ of the function $\tilde\Theta(\Lambda) := \Theta(\Lambda,a,0,(a-1)\Lambda)$ for some given initial value $a\in\R$. Then we have to solve the differential system \eqref{RedChar} with $u_1(0)=a$, $u_2(0)=0$, $\Lambda(0)=\Lambda_0$ on some interval $I$ with $0\in I$ by an appropriate numerical method, and we follow the path of this characteristic curve until we get $\Lambda(t_0)=0$ for some $t_0\in I$. Finally, we arrive at an eigenvalue $\lambda = (\mu+1)(\mu-2 u_2(t_0))$ of the angular spheroidal wave equation \eqref{ASWE} for $\gamma^2 = u_1(t_0)$. In the appendix, this method is illustrated by a concrete numerical example, namely for the values $\mu=1$ and $a=5$.

The differential system \eqref{GenChar} in Theorem \ref{CharEqu} may also be regarded as the characteristic equation of a PDE for a function $\Lambda(u_1,u_2)$ with only two variables $u_1$ and $u_2$. The question arises whether one can associate a more simple linear Hamiltonian system with the spheroidal wave equation ($\beta=0$), which then yields a PDE for the eigenvalues that is less complicated than \eqref{QPDE}. In the remainder of this paper we address this problem, where we restrict our consideration to the case of the prolate spheroidal wave equation \eqref{CSWE} with $\beta=0$ and $\gamma^2>0$. For this purpose, we fix $u_3=\Lambda(u_1+u_2^2-1)-2(\mu+1)u_2$ in \eqref{GSWE} resp. \eqref{SysA}, and we assume $u_1=\gamma^2>0$, in which case also $u_1+u_2^2>0$ holds.

First, let us introduce new parameters that make the coefficient matrix of the differential system \eqref{SysA} appear more symmetric, and these are
\begin{equation} \label{NewPar}
v_1 := \frac{u_2}{\sqrt{u_1+u_2^2}},\quad v_2 := \sqrt{u_1+u_2^2},\quad
\zeta := \frac{\Lambda(u_1+u_2^2)-2\alpha u_2}{\sqrt{u_1+u_2^2}}
\end{equation}
where $v_1\in(-1,1)$ and $v_2\in(0,\infty)$. Note that conversely
\begin{equation*}
u_1+u_2^2=v_2^2,\quad u_2 = v_1 v_2,\quad \Lambda=\frac{\zeta+2\alpha v_1}{v_2},\quad u_3+\Lambda = (\zeta-2\alpha v_1)v_2
\end{equation*}
By means of the transformation
\begin{equation} \label{Scale}
\hat y(z) = \begin{pmatrix} 1 & 0 \\[1ex] 0 & v_2 \end{pmatrix} y(z)
\end{equation}
the differential system \eqref{GSWE} resp. \eqref{SysA} is equivalent to the $2\times 2$ system
\begin{equation*}
\hat y'(z) = \left(
\frac{1}{z}\begin{pmatrix} -\alpha & \zeta-2\alpha v_1 \\[1ex] 0 & \alpha \end{pmatrix} 
+ \frac{1}{1-z}\begin{pmatrix} \alpha & 0 \\[1ex] -\zeta-2\alpha v_1 & -\alpha \end{pmatrix}
+ \begin{pmatrix} 2 v_1 v_2 & 2 v_2 \\[1ex] -2 v_2 & -2 v_1 v_2 \end{pmatrix}\right)\hat y(z) 
\end{equation*}
which can be written in the form $\hat\tau\hat y = \zeta\hat y$, where
\begin{equation*}
\hat\tau\hat y := W(z)^{-1}\left(J\hat y'(z) - \hat H(z)\hat y(z)\right),\quad
\hat H(z) := \begin{pmatrix} 2v_2+\frac{2\alpha v_1}{1-z} & -\frac{\alpha}{z}+\frac{\alpha}{1-z}+2v_1v_2 \\[1ex] 
-\frac{\alpha}{z}+\frac{\alpha}{1-z}+2v_1v_2 & 2v_2-\frac{2\alpha v_1}{z} \end{pmatrix}
\end{equation*}
By a similar reasoning as for \eqref{tau}, it follows that the differential expression $\hat\tau$ is in the limit point case at $z=0$ and $z=1$. Hence, the closure of its associated minimal operator is the only self-adjoint realization of $\hat\tau$ in the Hilbert space $\LzW$. We denote this self-adjoint operator by $\hat T=\hat T(v_1,v_2)$. By virtue of \eqref{Scale}, we have $\langle y,y\rangle_W < \infty$ if and only if $\langle\hat y,\hat y\rangle_W < \infty$. Hence, $\Lambda$ is an eigenvalue of $T(\U)$ if and only if $\zeta$ is an eigenvalue of $\hat T(v_1,v_2)$.

Now, like in Theorem \ref{Main}, suppose that $\zeta=\zeta(v_1,v_2)$ is an eigenvalue of $\hat T(v_1,v_2)$, which depends analytically on the parameters $(v_1,v_2)\in V$ in some domain $V\subset(-1,1)\times(0,\infty)$. If we introduce the matrix function
\begin{equation*}
\hat G(z,\zeta,v_1,v_2) := \begin{pmatrix} 
v_2 z & v_1 v_2 z + \alpha(v_1^2-1) - \frac{1}{2}v_1 v_2 \\[1ex] 
v_1 v_2 z + \alpha(v_1^2-1) - \frac{1}{2}v_1 v_2 & v_2(z-1) \end{pmatrix}
\end{equation*}
then a direct calculation shows that the deformation equation
\begin{equation*}
\pD{\hat G}{z} + (\zeta W + \hat H)J\hat G - \hat GJ(\zeta W + \hat H) 
= \hat f_1\,\pD{H}{v_1} + \hat f_2\,\pD{H}{v_2} + \hat g W
\end{equation*}
is fulfilled with the functions
\begin{gather*}
\hat f_1 := (\zeta+v_2)(1-v_1^2),\quad
\hat f_2 := (\zeta v_1+v_1v_2+\tfrac{1}{2})v_2,\quad
\hat g   := 4\alpha^2v_1(1-v_1^2)-\zeta v_1 v_2
\end{gather*}
From \cite[Theorem 2.2]{Schmid:2021} it follows that the eigenvalues $\zeta(v_1,v_2)$ satisfy the first-order quasilinear PDE
\begin{equation} \label{RPDE}
(\zeta+v_2)(1-v_1^2)\frac{\partial\zeta}{\partial v_1} 
+ (\zeta v_1+v_1v_2+\tfrac{1}{2})v_2\frac{\partial\zeta}{\partial v_2} 
= 4\alpha^2v_1(1-v_1^2)-\zeta v_1 v_2
\end{equation}
Using Lemma \ref{Embed} with $u_3=\Lambda(u_1+u_2^2-1)-2(\mu+1)u_2$, we obtain that $\lambda=\Lambda(u_1+u_2^2-1)+(\mu+1)(\mu-2u_2)$ is an eigenvalue of \eqref{CSWE} for the parameters $\gamma^2=u_1$, $\beta = -\Lambda(u_1+u_2^2-1)$ if and only if $\Lambda=0$ is an eigenvalue of $T(\U)$. In this case $\lambda=(\mu+1)(\mu-2u_2)$ and $\beta = 0$, i.e., the CSWE reduces to the angular spheroidal wave equation \eqref{ASWE}. Note that $\Lambda=0$ is equivalent to $\zeta = -2\alpha v_1 = -(\mu+1)v_1$. Since $u_1=(1-v_1^2)v_2^2$ and $u_2 = v_1 v_2$, it follows that $\lambda=(\mu+1)(\mu-2 v_1v_2)$ is an eigenvalue of \eqref{ASWE} for $\gamma^2=(1-v_1^2)v_2^2$ if and only if $\zeta = -(\mu+1)v_1$ is an eigenvalue of $\hat T(v_1,v_2)$.

Finally, let us simplify the PDE \eqref{RPDE} even further. To this end, we need to modify the parameters introduced in \eqref{NewPar} once more. We set
\begin{equation*}
v_1 = \tanh u,\quad v_2 = \E^t\cosh u,\quad\zeta = \omega-\E^t\cosh u
\end{equation*}
With these parameters $(t,u)\in\R^2$ we convert the eigenvalue problem $\hat\tau\hat y = \zeta\hat y$ to
\begin{equation} \label{HSWE}
W(z)^{-1}\left(J\hat y'(z) - \Phi(z,t,u)\hat y(z)\right) = \omega\hat y(z),\quad 0<z<1
\end{equation}
where now $\omega$ is the eigenvalue parameter, and
\begin{equation*}
\Phi(z,t,u) := 
\begin{pmatrix} 2\E^t\cosh u + \frac{2\alpha\tanh u - \E^t\cosh u}{1-z} & -\frac{\alpha}{z}+\frac{\alpha}{1-z}+2\E^t\sinh u \\[1ex] -\frac{\alpha}{z}+\frac{\alpha}{1-z}+2\E^t\sinh u & 2\E^t\cosh u - \frac{2\alpha\tanh u + \E^t\cosh u}{z} \end{pmatrix}
\end{equation*}

\begin{Theorem} \label{SWEPDE} 
Let $\omega=\omega(t,u)$ be an eigenvalue of the linear Hamiltonian system \eqref{HSWE}, and suppose that $\omega(t,u)$ depends analytically on the parameters $(t,u)\in S$ in some domain $S\subset\R^2$. Then $\omega(t,u)$ solves the quasilinear partial differential equation
\begin{equation} \label{SPDE}
\pD{\omega}{t} + 2\omega\cdot\pD{\omega}{u} 
= 2(\mu+1)^2\frac{\tanh u}{\cosh^2 u} + \E^{2t}\sinh(2u) + \E^t\cosh u
\end{equation}
on $S$. Moreover, $\lambda=(\mu+1)(\mu-2\E^t\sinh u)$ is an eigenvalue of the spheroidal wave equation \eqref{ASWE} for $\gamma=\E^t$ if and only if the equation $\omega(t,u)=\E^t\cosh u-(\mu+1)\tanh u$ is satisfied for some $(t,u)\in S$.
\end{Theorem}

\begin{Proof}
Note that $\omega = \zeta + v_2$, and that $v_1$ does not depend on $t$. Hence,
\begin{align*}
\pD{\omega}{t} 
& = \pD{\zeta}{v_2}\pD{v_2}{t}+\pD{v_2}{t}
  = \E^t\cosh u\,\pD{\zeta}{v_2} + \E^t\cosh u \\
\pD{\omega}{u} 
& = \pD{\zeta}{v_1}\pD{v_1}{u} + \pD{\zeta}{v_2}\pD{v_2}{u} + \pD{v_2}{u}
  = \frac{1}{\cosh^2 u}\,\pD{\zeta}{v_1} + \E^t\sinh u\,\pD{\zeta}{v_2}+\E^t\sinh u
\end{align*}
and therefore
\begin{align*}
\pD{\omega}{t} + 2\omega\cdot\pD{\omega}{u}
& = \frac{2\omega}{\cosh^2 u}\,\pD{\zeta}{v_1} 
  + (\E^t\cosh u + 2\omega\,\E^t\sinh u)\pD{\zeta}{v_2}
  + \E^t\cosh u + 2\omega\,\E^t\sinh u
\end{align*}
From
\begin{equation*}
\frac{2\omega}{\cosh^2 u} = 2(\zeta+v_2)(1-v_1^2),\quad
\E^t\cosh u + 2\omega\,\E^t\sinh u = 2(\zeta v_1+v_1v_2+\tfrac{1}{2})v_2
\end{equation*}
and the partial differential equation \eqref{RPDE} it follows that
\begin{align*}
\pD{\omega}{t} + 2\omega\cdot\pD{\omega}{u}
& = 2(\zeta+v_2)(1-v_1^2)\frac{\partial\zeta}{\partial v_1} 
  + 2(\zeta v_1+v_1v_2+\tfrac{1}{2})v_2\pD{\zeta}{v_2} + 2(\zeta v_1+v_1v_2+\tfrac{1}{2})v_2 \\
& = 8\alpha^2 v_1(1-v_1^2)-2\zeta v_1 v_2 + 2(\zeta v_1+v_1v_2+\tfrac{1}{2})v_2 \\
& = 8\alpha^2 v_1(1-v_1^2) + (2v_1 v_2+1)v_2 
  = 2(\mu+1)^2\frac{\tanh u}{\cosh^2 u} + 2\E^{2t}\sinh u\cosh u + \E^t\cosh u
\end{align*}
and this is exactly the PDE \eqref{SPDE}. As we have already noted, $\lambda=(\mu+1)(\mu-2 v_1v_2)=(\mu+1)(\mu-2\E^t\sinh u)$ is an eigenvalue of \eqref{ASWE} for $\gamma^2=(1-v_1^2)v_2^2=\E^{2t}$ if and only if $\omega=\zeta + v_2 = -(\mu+1)\tanh u+\E^t\cosh u$ is an eigenvalue of $\hat T(v_1,v_2)$ and hence of \eqref{HSWE}.
\end{Proof}

It is worth mentioning that the partial differential equation \eqref{SPDE} is basically a forced inviscid Burgers' equation. In fact, if we set $t=\frac{s}{2}$, then \eqref{SPDE} becomes
\begin{equation*}
\pD{\omega}{s} + \omega\cdot\pD{\omega}{u} 
= (\mu+1)^2\frac{\tanh u}{\cosh^2 u} + \tfrac{1}{2}\E^{s}\sinh(2u) + \tfrac{1}{2}\E^{\frac{s}{2}}\cosh u
\end{equation*}
and this quasilinear PDE can be written in the form $\pD{\omega}{s}+\omega\cdot\pD{\omega}{u} = f(s,u)$. The forced inviscid Burgers' equation is commonly encountered in fluid mechanics or gas dynamics, where it is used, for example, to study nonlinear interactions of dispersive waves. It is quite surprising that this PDE is related to the linear Hamiltonian system \eqref{HSWE} and thus as well to the eigenvalues of the spheroidal wave equation.

\section{Conclusion}
\label{sec:Conclusion}

Calculating the eigenvalues of a spheroidal wave equation is known to be a tricky task. There is an ongoing effort in finding new techniques to compute these eigenvalues in a convenient way and with high precision (see, for example, \cite{Falloon:2001}, \cite{Schmid:2023}, \cite{Flammer:1957}, \cite{Hodge:1970}, \cite{Kirby:2006}, \cite{Skoro:2015} and the references therein). In the present work, we used an approach via a $2\times 2$ differential system. Apart from some numerical issues, we have mainly studied the analytical dependence of the eigenvalues on the parameters $\gamma$ and $\beta$ (or $u_1$, $u_2$, $u_3$, respectively). Results such as those in Theorem \ref{SWEPDE} and relations like the PDE \eqref{QPDE} may help to obtain new insights into the analytic properties of the angular or Coulomb spheroidal eigenvalues. The CSWE, in turn, is a special case of the confluent Heun differential equation (CHE), which is likewise a subject of current research. The methods presented here may also apply to the CHE and to even more general eigenvalue problems.

\section*{Appendix: A numerical example}

\begin{figure}
\centering
\includegraphics{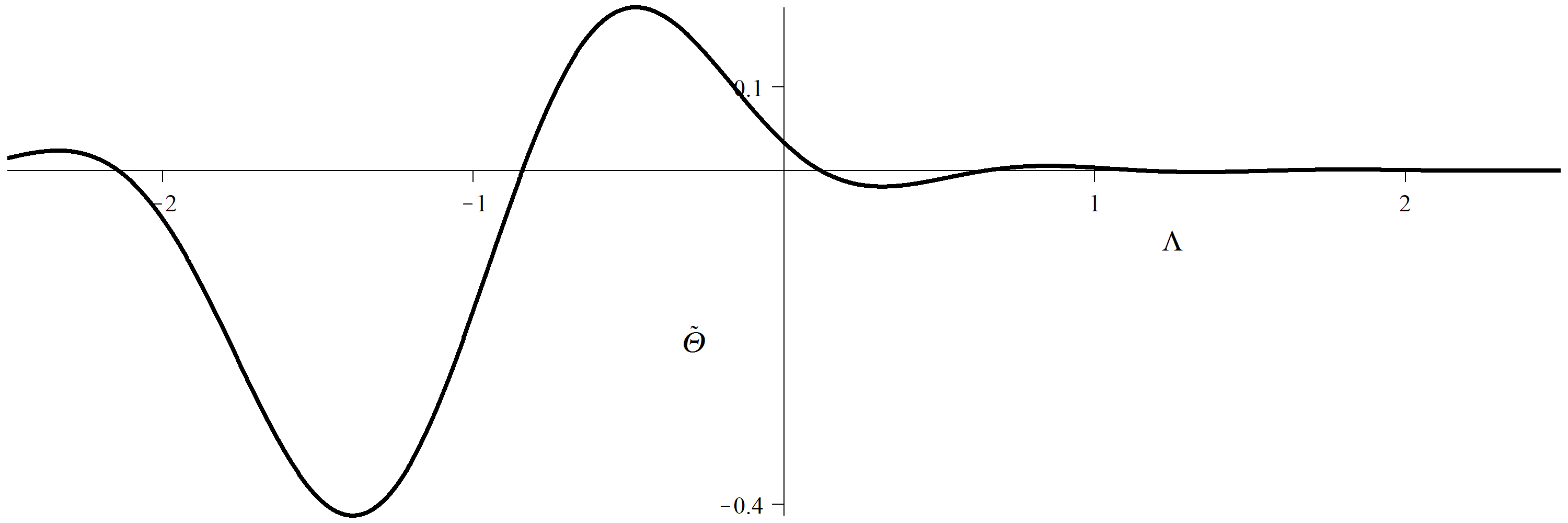}
\caption{The function $\tilde\Theta(\Lambda)$ for $\mu=1$ and $a=5$ with a zero at $\Lambda_0=-0.841720\ldots$} \label{fig:Theta}
\end{figure}

\begin{figure}
\centering
\includegraphics{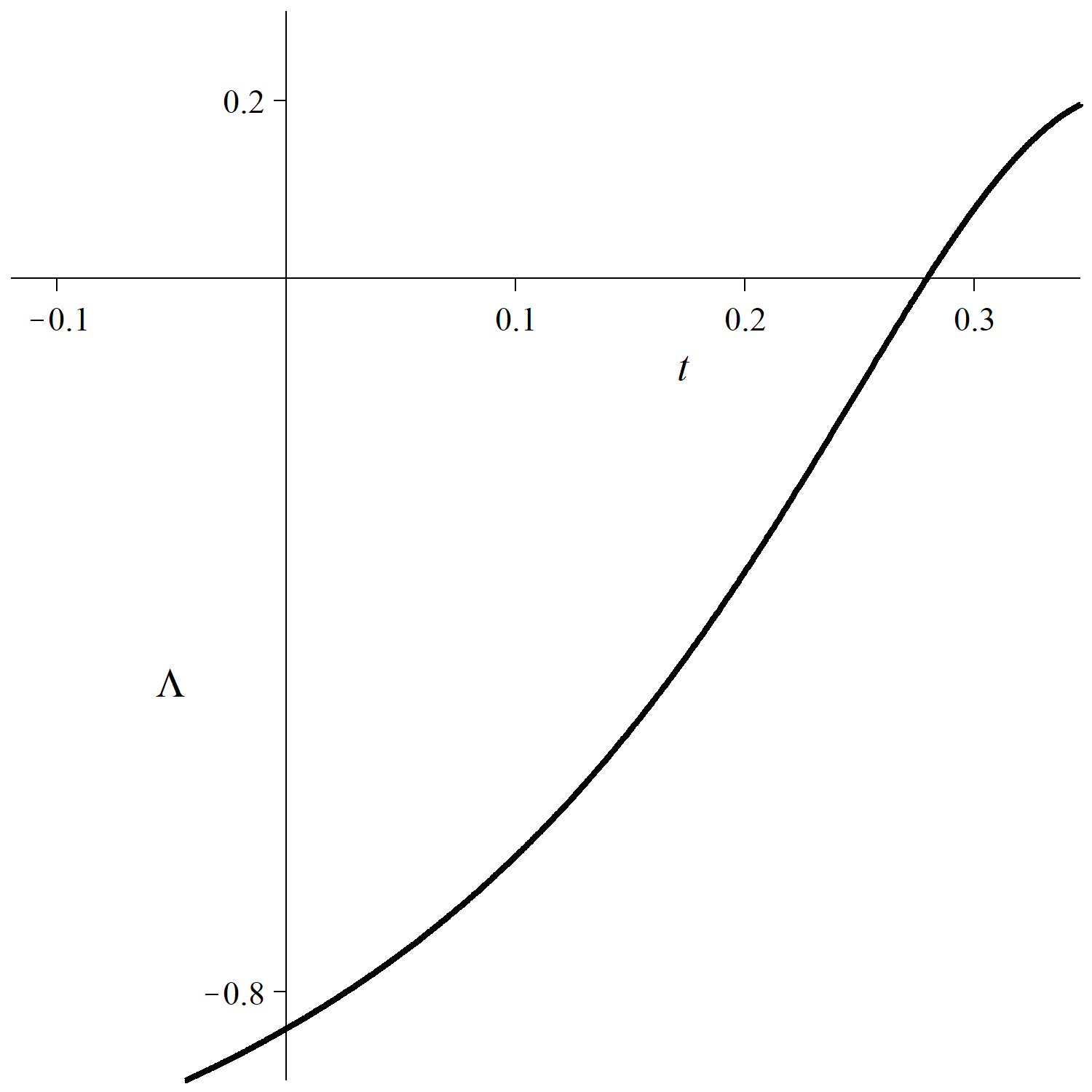}
\caption{The solution $\Lambda(t)$ of \eqref{RedChar} for $\mu=1$ starting from $u_1(0)=5$, $u_2(0)=0$ and $\Lambda(0)=-0.841720\ldots$} \label{fig:Lambda}
\end{figure}

In section \ref{sec:Evaluation}, as an application of Theorem \ref{CharEqu}, it has been suggested to compute an eigenvalue of the spheroidal wave equation by starting from a zero of the function $\tilde\Theta$ and then following a characteristic curve of the PDE \eqref{QPDE} until a zero of the function $\Lambda(t)$ is located. As a concrete numerical example, let us consider the function $\tilde\Theta(\Lambda) := \Theta(\Lambda,a,0,(a-1)\Lambda)$ in the case $\mu=1$ for $a=5$. The limit formula in Lemma \ref{Limit} provides the function shown in Fig. \ref{fig:Theta}. It has a zero at $\Lambda_0 =-0.8417200168449013$, which can be determined, for instance, with the secant method. Here, all results are displayed rounded to $16$ decimal places, while the calculations were internally done with higher precision. Solving the differential system \eqref{RedChar} for the initial values $u_1(0)=5$, $u_2(0)=0$, $\Lambda(0)=\Lambda_0$ by means of the Runge-Kutta method, we obtain $u_1(t)$, $u_2(t)$ and $\Lambda(t)$. The function $\Lambda(t)$ is illustrated in Fig. \ref{fig:Lambda}; it has a zero at $t_0=0.2793371978706399$. Evaluating  $u_1$, $u_2$ at $t_0$ finally gives the eigenvalue $\lambda=2(1-2 u_2(t_0)) =-5.2736106330552739$ for the spheroidal wave equation \eqref{ASWE} with parameter $\gamma^2 = u_1(t_0) = 8.7417666942941543$.

\end{document}